\documentclass[12pt]{article}

\usepackage{amsmath,amssymb,amsthm,amscd,a4wide,upref}
\usepackage[toc]{appendix}

\usepackage{hyperref}
\hypersetup{colorlinks,citecolor=blue,filecolor=black,linkcolor=blue,urlcolor=blue}
\usepackage{enumerate}
\usepackage{mathtools}

\begin{document}

\newcommand{\ad}{{\rm ad}}
\newcommand{\cri}{{\rm cri}}
\newcommand{\End}{{\rm{End}\ts}}
\newcommand{\Rep}{{\rm{Rep}\ts}}
\newcommand{\Hom}{{\rm{Hom}}}
\newcommand{\Mat}{{\rm{Mat}}}
\newcommand{\ch}{{\rm{ch}\ts}}
\newcommand{\chara}{{\rm{char}\ts}}
\newcommand{\diag}{{\rm diag}}
\newcommand{\non}{\nonumber}
\newcommand{\wt}{\widetilde}
\newcommand{\wh}{\widehat}
\newcommand{\ot}{\otimes}
\newcommand{\la}{\lambda}
\newcommand{\La}{\Lambda}
\newcommand{\De}{\Delta}
\newcommand{\al}{\alpha}
\newcommand{\be}{\beta}
\newcommand{\ga}{\gamma}
\newcommand{\Ga}{\Gamma}
\newcommand{\ep}{\epsilon}
\newcommand{\ka}{\kappa}
\newcommand{\vk}{\varkappa}
\newcommand{\si}{\sigma}
\newcommand{\vs}{\varsigma}
\newcommand{\vp}{\varphi}
\newcommand{\de}{\delta}
\newcommand{\ze}{\zeta}
\newcommand{\om}{\omega}
\newcommand{\Om}{\Omega}
\newcommand{\ee}{\epsilon^{}}
\newcommand{\su}{s^{}}
\newcommand{\hra}{\hookrightarrow}
\newcommand{\ve}{\varepsilon}
\newcommand{\ts}{\,}
\newcommand{\vac}{\mathbf{1}}
\newcommand{\di}{\partial}
\newcommand{\qin}{q^{-1}}
\newcommand{\tss}{\hspace{1pt}}
\newcommand{\Sr}{ {\rm S}}
\newcommand{\U}{ {\rm U}}
\newcommand{\BL}{ {\overline L}}
\newcommand{\BE}{ {\overline E}}
\newcommand{\BP}{ {\overline P}}
\newcommand{\AAb}{\mathbb{A}\tss}
\newcommand{\CC}{\mathbb{C}\tss}
\newcommand{\KK}{\mathbb{K}\tss}
\newcommand{\QQ}{\mathbb{Q}\tss}
\newcommand{\SSb}{\mathbb{S}\tss}
\newcommand{\TT}{\mathbb{T}\tss}
\newcommand{\ZZ}{\mathbb{Z}\tss}
\newcommand{\DY}{ {\rm DY}}
\newcommand{\X}{ {\rm X}}
\newcommand{\Y}{ {\rm Y}}
\newcommand{\Z}{{\rm Z}}
\newcommand{\Ac}{\mathcal{A}}
\newcommand{\Lc}{\mathcal{L}}
\newcommand{\Mc}{\mathcal{M}}
\newcommand{\Pc}{\mathcal{P}}
\newcommand{\Qc}{\mathcal{Q}}
\newcommand{\Rc}{\mathcal{R}}
\newcommand{\Sc}{\mathcal{S}}
\newcommand{\Tc}{\mathcal{T}}
\newcommand{\Bc}{\mathcal{B}}
\newcommand{\Ec}{\mathcal{E}}
\newcommand{\Fc}{\mathcal{F}}
\newcommand{\Gc}{\mathcal{G}}
\newcommand{\Hc}{\mathcal{H}}
\newcommand{\Uc}{\mathcal{U}}
\newcommand{\Vc}{\mathcal{V}}
\newcommand{\Wc}{\mathcal{W}}
\newcommand{\Yc}{\mathcal{Y}}
\newcommand{\Ar}{{\rm A}}
\newcommand{\Br}{{\rm B}}
\newcommand{\Ir}{{\rm I}}
\newcommand{\Fr}{{\rm F}}
\newcommand{\Jr}{{\rm J}}
\newcommand{\Or}{{\rm O}}
\newcommand{\GL}{{\rm GL}}
\newcommand{\Spr}{{\rm Sp}}
\newcommand{\Rr}{{\rm R}}
\newcommand{\Zr}{{\rm Z}}
\newcommand{\gl}{\mathfrak{gl}}
\newcommand{\middd}{{\rm mid}}
\newcommand{\ev}{{\rm ev}}
\newcommand{\Pf}{{\rm Pf}}
\newcommand{\Norm}{{\rm Norm\tss}}
\newcommand{\oa}{\mathfrak{o}}
\newcommand{\spa}{\mathfrak{sp}}
\newcommand{\osp}{\mathfrak{osp}}
\newcommand{\f}{\mathfrak{f}}
\newcommand{\g}{\mathfrak{g}}
\newcommand{\h}{\mathfrak h}
\newcommand{\n}{\mathfrak n}
\newcommand{\z}{\mathfrak{z}}
\newcommand{\Zgot}{\mathfrak{Z}}
\newcommand{\p}{\mathfrak{p}}
\newcommand{\sll}{\mathfrak{sl}}
\newcommand{\agot}{\mathfrak{a}}
\newcommand{\qdet}{ {\rm qdet}\ts}
\newcommand{\Ber}{ {\rm Ber}\ts}
\newcommand{\HC}{ {\mathcal HC}}
\newcommand{\cdet}{{\rm cdet}}
\newcommand{\rdet}{{\rm rdet}}
\newcommand{\tr}{ {\rm tr}}
\newcommand{\gr}{ {\rm gr}\ts}
\newcommand{\str}{ {\rm str}}
\newcommand{\loc}{{\rm loc}}
\newcommand{\Gr}{{\rm G}}
\newcommand{\sgn}{ {\rm sgn}\ts}
\newcommand{\sign}{{\rm sgn}}
\newcommand{\ba}{\bar{a}}
\newcommand{\bb}{\bar{b}}
\newcommand{\bi}{\bar{\imath}}
\newcommand{\bj}{\bar{\jmath}}
\newcommand{\bk}{\bar{k}}
\newcommand{\bl}{\bar{l}}
\newcommand{\hb}{\mathbf{h}}
\newcommand{\Sym}{\mathfrak S}
\newcommand{\fand}{\quad\text{and}\quad}
\newcommand{\Fand}{\qquad\text{and}\qquad}
\newcommand{\For}{\qquad\text{or}\qquad}
\newcommand{\OR}{\qquad\text{or}\qquad}
\newcommand{\grpr}{{\rm gr}^{\tss\prime}\ts}
\newcommand{\degpr}{{\rm deg}^{\tss\prime}\tss}

\numberwithin{equation}{section}

\newtheorem{thm}{Theorem}[section]
\newtheorem{lem}[thm]{Lemma}
\newtheorem{prop}[thm]{Proposition}
\newtheorem{cor}[thm]{Corollary}
\newtheorem{conj}[thm]{Conjecture}
\newtheorem*{mthm}{Main Theorem}
\newtheorem*{mthma}{Theorem A}
\newtheorem*{mthmb}{Theorem B}
\newtheorem*{mthmc}{Theorem C}
\newtheorem*{mthmd}{Theorem D}

\theoremstyle{definition}
\newtheorem{defin}[thm]{Definition}

\theoremstyle{remark}
\newtheorem{remark}[thm]{Remark}
\newtheorem{example}[thm]{Example}

\newcommand{\bth}{\begin{thm}}
\renewcommand{\eth}{\end{thm}}
\newcommand{\bpr}{\begin{prop}}
\newcommand{\epr}{\end{prop}}
\newcommand{\ble}{\begin{lem}}
\newcommand{\ele}{\end{lem}}
\newcommand{\bco}{\begin{cor}}
\newcommand{\eco}{\end{cor}}
\newcommand{\bde}{\begin{defin}}
\newcommand{\ede}{\end{defin}}
\newcommand{\bex}{\begin{example}}
\newcommand{\eex}{\end{example}}
\newcommand{\bre}{\begin{remark}}
\newcommand{\ere}{\end{remark}}
\newcommand{\bcj}{\begin{conj}}
\newcommand{\ecj}{\end{conj}}

\newcommand{\bal}{\begin{aligned}}
\newcommand{\eal}{\end{aligned}}
\newcommand{\beq}{\begin{equation}}
\newcommand{\eeq}{\end{equation}}
\newcommand{\ben}{\begin{equation*}}
\newcommand{\een}{\end{equation*}}

\newcommand{\bpf}{\begin{proof}}
\newcommand{\epf}{\end{proof}}

\def\beql#1{\begin{equation}\label{#1}}


\newcommand{\Res}{\mathop{\mathrm{Res}}}

\title{\Large\bf Center of the quantum affine vertex algebra\\[0.2em]
associated with trigonometric
$R$-matrix}

\author{{Slaven Ko\v{z}i\'{c}\quad and\quad Alexander Molev}}

\date{} 
\maketitle

\vspace{5 mm}

\begin{abstract}
We consider the quantum vertex algebra associated with the trigonometric $R$-matrix
in type $A$ as defined by Etingof and Kazhdan. We show that its
center is a commutative associative algebra and construct
an algebraically independent family of topological generators of the center
at the critical level.



\end{abstract}

\vspace{5 mm}

%

\section{Introduction}
\label{sec:int}

A general definition of {\em quantum vertex algebra} was given by
P.~Etingof and D.~Kazhdan~\cite{ek:ql5}. In particular,
a {\em quantum affine vertex algebra} can be associated with a rational,
trigonometric or elliptic $R$-matrix.
A suitably normalized Yang $R$-matrix gives rise to a
quantum vertex algebra structure on the vacuum module $\Vc_c(\gl_n)$
over the double Yangian
$\DY(\gl_n)$. In our previous paper~\cite{jkmy:cq} coauthored with N.~Jing and F.~Yang,
we introduced the {\em center} of an arbitrary quantum vertex algebra and
described the center $\z\big(\Vc_c(\gl_n)\big)$
of the quantum affine vertex algebra $\Vc_c(\gl_n)$.
We showed that the center at the critical level $c=-n$
possesses large families of algebraically independent topological generators
in a complete analogy with the affine vertex algebra \cite{ff:ak}; see also
\cite{f:lc}.

Our goal in this paper is to give a similar description of the center of
the quantized universal enveloping algebra $\U(R)$, associated with a normalized
trigonometric $R$-matrix $R$, as a quantum vertex algebra.
We show that the center
of the $h$-adically completed quantum affine vertex algebra $\U_c(R)$
at the level $c\in\CC$ is a commutative algebra. Moreover, we produce
a family of algebraically independent
topological generators of the center in an explicit form.
We show
that taking their `classical limits' reproduces the corresponding generators of
the center of
the quantum affine vertex algebra $\Vc_c(\gl_n)$.
In particular,
as with the rational case, the center of the $h$-adic completion of $\U_c(R)$
is `large' at the critical level $c=-n$, and trivial otherwise.

Despite an apparent analogy between the rational and trigonometric cases,
there are significant differences in the constructions. In the rational case,
the quantum vertex algebra structure is essentially determined by that of the double
Yangian $\DY(\gl_n)$. One could expect that the role of the double
Yangian in the trigonometric case to be played by the quantum affine algebra $\U_q(\gl_n)$.
In fact, as explained in \cite{ek:ql5},
a more subtle structure is to be used instead.
Nonetheless,
the technical part is quite similar to that of the paper \cite{fjmr:hs}, where
explicit constructions of elements of the center
of the completed quantum affine algebra were given, and which stems from
the pioneering work of
N.~Reshetikhin and M.~Semenov-Tian-Shansky~\cite{rs:ce}; see also
J.~Ding and P.~Etingof~\cite{de:cq} and
E.~Frenkel and N.~Reshetikhin~\cite{fr:qa}.

\section{Quantized universal enveloping algebra}
\label{sec:QUE}

In accordance to \cite{ek:ql5}, a normalized $R$-matrix is needed
to define an appropriate version
of the quantized universal enveloping algebra $\U(R)$.
Namely, the $R$-matrix should satisfy the {\em unitarity}
and {\em crossing symmetry} properties.
Its existence is established in
\cite[Proposition~1.2]{ek:ql4}. The normalizing factor does not admit
a simple closed expression. A description of this factor in the rational
case is also given in \cite[Section~2.2]{jkmy:cq}. Here we give a similar
description in the trigonometric case which also implies the existence
of the normalized $R$-matrix.  We start by recalling some standard tensor notation.

We let $e_{ij}\in\End\CC^n$ denote the standard matrix units.
For an element
\ben
C=\sum_{i,j,r,s=1}^n c^{}_{ijrs}\ts e_{ij}\ot e_{rs}\in
\End \CC^n\ot\End \CC^n,
\een
and any two indices $a,b\in\{1,\dots,m\}$ such that $a\ne b$,
we denote by $C_{a\tss b}$ the element of the algebra $(\End\CC^n)^{\ot m}$ with $m\geqslant 2$
given by
\beql{cab}
C_{a\tss b}=\sum_{i,j,r,s=1}^n c^{}_{ijrs}\ts (e_{ij})_a\tss (e_{rs})_b,
\qquad
(e_{ij})_a=1^{\ot(a-1)}\ot e_{ij}\ot 1^{\ot(m-a)}.
\eeq
We regard the matrix transposition as the linear map
\ben
t:\End\CC^n\to\End\CC^n,\qquad e_{ij}\mapsto e_{ji}.
\een
For any $a\in\{1,\dots,m\}$ we will denote by $t_a$ the corresponding
partial transposition on the algebra $(\End\CC^n)^{\ot m}$ which acts as $t$ on the
$a$-th copy of $\End \CC^n$ and as the identity map on all the other tensor factors.

Introduce the two-parameter $R$-matrix
$\overline{R}(u,v)\in\End\CC^n\otimes\End\CC^n  [[u,v,h]]$ as a formal power series
in the variables $u,v,h$ by
\begin{align}
\overline{R}(u,v)=&\big(e^{u-h/2}-e^{v+h/2}\big)\sum_{i} e_{ii}\otimes e_{ii}
+\big(e^u -e^v\big)\sum_{i\neq j} e_{ii}\otimes e_{jj}\nonumber \\
&+ \big(e^{-h/2}-e^{h/2}\big)\tss e^u \sum_{i>j} e_{ij}\otimes e_{ji}
+ \big(e^{-h/2}-e^{h/2}\big)\tss e^v\sum_{i<j} e_{ij}\otimes e_{ji},\label{R2p}
\end{align}
where the summation indices run over the set $\{1,\dots,n\}$.
We will also use the one-parameter $R$-matrix
$\overline{R}(u)\in\End\CC^n\otimes\End\CC^n \ts((u))[[h]]$
defined by
\begin{align}
\overline{R}(u)=\frac{\overline{R}(u,0)}{e^{u-h/2}-e^{h/2}}=&\sum_{i} e_{ii}\otimes e_{ii}
+e^{-h/2}\ts\frac{1-e^u}{1-e^{u-h}}\tss\sum_{i\neq j} e_{ii}\otimes e_{jj}\nonumber \\
&+ \frac{\big(1-e^{-h}\big)\ts e^u}{1-e^{u-h}}\sum_{i>j} e_{ij}\otimes e_{ji}
+ \frac{1-e^{-h}}{1-e^{u-h}}\sum_{i<j} e_{ij}\otimes e_{ji}.\label{R1p}
\end{align}
Here and below expressions of the form $\big(1-e^{u+a\tss h}\big)^{-1}$ with $a\in\CC$
should be understood as elements of the algebra $\CC((u))[[h]]$,
\ben
\big(1-e^{u+a\tss h}\big)^{-1}
=-u^{-1}\Bigg(\sum_{l\geqslant 1}\frac{(u+a\tss h)^{l-1}}{l!}\Bigg)^{-1}\big(1+a\tss h/u\big)^{-1}
\ts\in\ts\CC((u))[[h]].
\een
Denote by $D$  the diagonal $n \times n$ matrix
\beql{diagonal}
D = \diag\ts\big[e^{\frac{(n-1)\tss h}{2}},e^{\frac{(n-3)\tss h}{2}},\dots,e^{-\frac{(n-1)\tss h}{2}}\big]
\eeq
with entries in $\CC[[h]]$.
The following proposition goes back to \cite[Proposition~4.7]{fri:qa}
and is a version of \cite[Proposition~1.2]{ek:ql4}. We use the notation \eqref{cab}.

\bpr\label{prop:normal}
There exists a unique series $g(u)\in 1+h\tss\CC((u))[[h]]$ $(\text{depending on\ } $n$)$
such that the $R$-matrix $R(u)=g(u)\tss\overline{R}(u)$ possesses  the unitarity property
\beql{uni}
R_{12}(u) R_{21}(-u) =1
\eeq
and the crossing symmetry properties
\beql{csym}
\big(R_{12}(u)^{-1}\big)^{t_2}D_2 R_{12}(u+nh)^{t_2}=D_2\fand
R_{12}(u+nh)^{t_1}D_1 \big( R_{12}(u)^{-1}\big)^{t_1}=D_1.
\eeq
\epr

\bpf
Due to the well-known properties of the $R$-matrix \eqref{R1p} (see \cite{fri:qa}),
identities \eqref{uni} and \eqref{csym} will hold for
$R(u)=g(u)\overline{R}(u)$
if and only if $g(u)$ satisfies the relations
\beql{g6}
g(u)\tss g(-u)=1
\eeq
and
\beql{g1}
g(u+nh)=g(u)\ts\frac{\big(1-e^{u+h}\big)\big(1-e^{u+(n-1)h}\big)}{\big(1-e^u\big)\big(1-e^{u+nh}\big)}.
\eeq
It is well known by \cite{fri:qa} that
there exists a unique formal power series $f(x)\in\CC(q)[[x]]$ of the form
\beql{fxdef}
f(x)=1+\sum_{k=1}^{\infty}f_kx^k,\qquad f_k=f_k(q),
\eeq
whose coefficients $f_k$ are determined
by the relation
\beql{fx}
f(xq^{2n})=f(x)\ts\frac{(1-xq^2)\tss(1-xq^{2n-2})}{(1-x)\tss(1-xq^{2n})}.
\eeq
Equivalently, the series \eqref{fxdef} is a unique solution of the equation
\beql{fxeq}
f(x)f(xq^2)\dots f(xq^{2n-2})=\frac{1-x}{1-x\tss q^{2n-2}}.
\eeq
Observe that $f(x)$ admits the presentation
\beql{fxax}
f(x)=1+\sum_{k=1}^{\infty}a_k\Big(\frac{x}{1-x}\Big)^k,
\eeq
where all rational functions $a_k/(q-1)^k\in\CC(q)$ are regular at $q=1$. Indeed, this claim
means that all coefficients of the formal series
\ben
b(z)=1+\sum_{k=1}^{\infty}\frac{a_k}{(q-1)^k}\tss z^k
\een
are regular at $q=1$. However, rewriting the equation \eqref{fxeq} in terms of $b(z)$ we get
\ben
b(z)\ts b\Big(\frac{zq^2}{1-z\tss\frac{1-q^2}{1-q}}\Big)\dots
b\Big(\frac{zq^{2n-2}}{1-z\tss\frac{1-q^{2n-2}}{1-q}}\Big)
=\frac{1}{1-z\tss\frac{1-q^{2n-2}}{1-q}}.
\een
This implies a system of recurrence relations for the coefficients of the series $b(z)$ so that an easy
induction argument
shows that each of the coefficients is regular at $q=1$. Thus, making the substitution
\beql{xq}
x=e^u \Fand q=e^{h/2}
\eeq
in \eqref{fxax} we obtain a well-defined element $\wt g(u)\in 1+h\tss\CC((u))[[h]]$
satisfying \eqref{g1}.

Now set
\ben
\vp(u) = \wt g(u)\tss \wt g(-u) \in 1+h\tss\CC((u))[[h]].
\een
Since $\vp(u)=\vp(-u)$, the Laurent series
\ben
\vp(u)=\sum_{s\in\mathbb{Z}}\vp_s u^s,\qquad \vp_s\in\CC[[h]],
\een
contains only even powers of $u$. Moreover,
\eqref{g1} implies the relation
\beql{g9}
\vp(u)=\vp(u+nh).
\eeq
By comparing the coefficients of the negative powers of $u$ on both sides, we
conclude that $\vp_s=0$ for all $s<0$. Similarly,
by considering nonnegative powers of $u$ in \eqref{g9}
we get $\vp_s=0$ for all $s>0$ so that
$\vp(u)=\vp_0$ is an element of $1+h\tss\CC[[h]]$.
Let $\psi\in 1+h\tss\CC[[h]]$ be such that $\psi^2\tss \vp_0=1$. Then
the series $g(u)=\psi\tss \wt g(u)$ satisfies both \eqref{g6} and \eqref{g1}.

A direct argument with formal series shows that the conditions \eqref{g6} and \eqref{g1}
uniquely determine $g(u)$. The details are given in Appendix~\ref{app:gu} which also
contains a direct proof of the existence of $g(u)$.
\epf

The first few terms of the series $g(u)$ are found by
\ben
g(u)=1+\frac{\big(n-1\big)\big(1+e^u \big)}{2n \big(1-e^u \big)}h+\frac{\big(n-1\big)^2
\big(1+e^u \big)^2 }{8n^2 \big(1-e^u \big)^2} h^2+\dots.
\een

\bco\label{cor:gu}
The series $g(u)$ satisfies the relation
\beql{Geq}
g(u)\tss g(u+h)\dots g\big(u+(n-1)h\big)=e^{(n-1)h/2}\ts\frac{1-e^u}{1-e^{u+(n-1)h}}.
\eeq
\eco

\bpf
Denote by $G(u)$ the series on the left hand side of \eqref{Geq}. This is an
element of $1+h\tss\CC((u))[[h]]$ which satisfies
\ben
G(u)\tss G\big({-}u-(n-1)h\big)=1
\een
and
\ben
G(u+h)=G(u)\tss \frac{(1-e^{u+h})\big(1-e^{u+(n-1)h}\big)}{(1-e^u)(1-e^{u+nh})}
\een
by \eqref{g6} and \eqref{g1}, respectively. The same argument as in Appendix~\ref{app:gu}
shows that any series $G(u)\in 1+h\tss\CC((u))[[h]]$
satisfying these two properties is determined
uniquely. The series
on the right hand side of \eqref{Geq} also
satisfies these two properties and so
the claim follows.
\epf

\bre\label{rem:fxgu}
The $R$-matrices given in \eqref{R1p} and
in \cite[eq. (2.2)]{fjmr:hs} are related by the change of parameters
\eqref{xq}. Note, however, that although the relations
\eqref{g1} and \eqref{fx} correspond to each other under this change,
the difference equation \eqref{fx} determines $f(x)$ uniquely,
as an element of $1+x\tss\CC(q)[[x]]$, whereas $g(u)$ is regarded
as an element of a different algebra of power series, namely,
$1+h\tss\CC((u))[[h]]$. Moreover, as shown in Appendix~\ref{app:gu},
the property \eqref{g6} is necessary to
guarantee that $g(u)$ is determined uniquely.
\qed
\ere

To make a connection with the quantum affine vertex algebra associated with the rational
$R$-matrix (see Proposition~\ref{prop:isomgraded} below),
define the $\ZZ$-gradation on the algebra $\End\CC^n \ot \End\CC^n [u^{\pm 1},h]$ by
setting
\beql{filtr}
\deg u^k h^l=-k-l
\eeq
and assigning the zero degree to elements of $\End\CC^n \ot \End\CC^n$.
Extend the degree function \eqref{filtr}
to the algebra
of formal series $\End\CC^n \ot \End\CC^n ((u))[[h]]$
by allowing it to take the infinite value. Elements of finite degree
will then form a subalgebra and we denote it by
$\End\CC^n \ot \End\CC^n ((u))[[h]]_{\text{fin}}$.

The $R$-matrix $\overline{R}(u)$ defined in \eqref{R1p} belongs to
$\End\CC^n \ot \End\CC^n ((u))[[h]]_{\text{fin}}$ and its degree is zero.
Denote by $\overline{R}^{\ts\textrm{rat}}(u)$ its component of degree zero.
It is easy to check that  $\overline{R}^{\ts\textrm{rat}}(u)$ coincides with the
Yang $R$-matrix, up to a scalar factor:
\ben
\overline{R}^{\ts\textrm{rat}}(u)=\frac{u}{u-h} \Big(1-\frac{h}{u}P    \Big)
\ts\in\ts \End\CC^n \ot \End\CC^n [[h/u]],
\een
where $P\in\End\CC^n\ot\End\CC^n$ is the permutation operator given by
\beql{permut}
P = \sum_{i,j=1}^n e_{ij}\ot e_{ji}.
\eeq
Similarly,
denote the highest degree component of the $R$-matrix $R(u)$,
defined in Proposition~\ref{prop:normal}, by $R^{\tss\textrm{rat}}(u)$.
The unitarity property \eqref{uni}
 and the crossing symmetry properties \eqref{csym} for the $R$-matrix $R(u)$ imply
 \beql{unirat}
 R^{\tss\textrm{rat}}_{12}(u) R^{\tss\textrm{rat}}_{21}(-u) =1
 \eeq
 and
 \ben
 \left(R^{\tss\textrm{rat}}_{12}(u)^{-1}\right)^{t_2} R^{\tss\textrm{rat}}_{12}(u+nh)^{t_2}=1\Fand
 R^{\tss\textrm{rat}}_{12}(u+nh)^{t_1} \left( R^{\tss\textrm{rat}}_{12}(u)^{-1}\right)^{t_1}=1.
 \een
 Furthermore, let $g^{\tss\textrm{rat}}(u)$ be
the highest degree component of the series $g(u)$,
 defined in Proposition \ref{prop:normal}.
Then $R^{\ts\textrm{rat}}(u)=g^{\tss\textrm{rat}}(u)\overline{R}^{\ts\textrm{rat}}(u)$.
It is clear from the proof of Proposition~\ref{prop:normal} that
 $g^{\tss\textrm{rat}}(u)\in 1+(h/u)\tss\CC[[h/u]]$.
By taking the highest degree components in \eqref{g1} we get
 \beql{ga1}
g^{\tss\textrm{rat}}(u+nh) =g^{\tss\textrm{rat}}(u)\tss \frac{(u+h)\big(u+(n-1)h\big)}{u(u+nh)}.
 \eeq
 Since $1-h/u$ is invertible in $\CC[[h/u]]$, by  replacing
 $g^{\tss\textrm{rat}}(u)$ with $(1-h/u)\tss \overline{g}(u)$ in \eqref{ga1}
 we obtain the equivalent equation
 \beql{ga2}
 \overline{g}(u+nh)=\Big(1-\frac{h^2}{u^2} \Big)\ts\overline{g}(u).
 \eeq
 Equation \eqref{ga2} has a unique
 solution for the series $\overline{g}(u)$ in $ 1+(h/u)\CC[[h/u]]$;
 see e.g.~\cite[Section~2.2]{jkmy:cq}. Therefore,
\ben
 R^{\ts\textrm{rat}}(u)=g^{\tss\textrm{rat}}(u)
 \overline{R}^{\ts\textrm{rat}}(u)=\overline{g}(u)\Big(1-\frac{h}{u}P \Big).
\een

The {\em quantized universal enveloping algebra} $\U(R)$
is the associative algebra over $\CC[[h]]$ generated by  elements
$l_{ij}^{(-r)}$, where $1\leqslant i,j\leqslant n$ and $r=1,2,\dots$,
subject to the defining relations
\beql{RTT}
R(u-v)\tss L_{1}^{+}(u)\tss L_{2}^{+}(v)=L_{2}^{+}(v)\tss L_{1}^{+}(u)\tss R(u-v),
\eeq
where the matrix $L^{+}(u)$ is given by
\beql{T+}
L^+ (u)=\sum_{i,j=1}^n e_{ij}\otimes l_{ij}^{+}(u)\ts\in\ts \End\CC^n \otimes \U(R)[[u]]
\eeq
and
\ben
l_{ij}^{+}(u)=\delta_{ij}-h\tss\sum_{r=1}^\infty l_{ij}^{(-r)}u^{r-1}\ts\in\ts \U(R)[[u]].
\een
Here we extend the notation \eqref{cab}
to matrices of the
form \eqref{T+}.
A subscript indicates its copy in the multiple
tensor product algebra
\ben
\underbrace{\End\CC^n\ot\dots\ot\End\CC^n}_m\ot \U(R)[[u]],
\een
so that
\ben
L_a^+(u)=\sum_{i,j=1}^n 1^{\ot (a-1)}\ot e_{ij}\ot 1^{\ot (m-a)}\ot l_{ij}^+(u).
\een
We take $m=2$ for the defining relations \eqref{RTT}. Note that
the $R$-matrix $R(u)$ in \eqref{RTT} can be replaced with $\overline R(u)$
to define the same algebra
$\U(R)$.

Recall that
the {\em dual Yangian} $\Y^+(\gl_n)$ for $\gl_n$ is the associative
algebra over $\CC[[h]]$ generated by  elements
${t}_{ij}^{(-r)}$, where $1\leqslant i,j\leqslant n$ and $r=1,2,\dots$,
subject to the defining relations
\beql{RTTy}
R^{\ts\textrm{rat}}(u-v)\tss{T}_{1}^{+}(u)\tss{T}_{2}^{+}(v)
={T}_{2}^{+}(v)\tss{T}_{1}^{+}(u)\tss R^{\ts\textrm{rat}}(u-v),
\eeq
where the matrix ${T}^{+}(u)$ is given by
\ben
{T}^+ (u)=\sum_{i,j=1}^n e_{ij}\otimes
{t}_{ij}^{+}(u)\ts\in\ts \End\CC^n \otimes \Y^+(\gl_n)[[u]]
\een
and
\ben
{t}_{ij}^{+}(u)=\delta_{ij}-h\sum_{r=1}^\infty
{t}_{ij}^{(-r)}u^{r-1}\ts\in\ts \Y^+(\gl_n)[[u]].
\een

Introduce the descending filtration
\ben
{}\dots\supset \U^{(2)}\supset \U^{(1)}\supset \U^{(0)}\supset \U^{(-1)}\supset \U^{(-2)}\supset\cdots
\een
on
$\U(R)$ by setting
\beql{filtr2}
\deg h=-1\Fand \deg l_{ij}^{\tss(-r)}=r
\eeq
so that for any $r\in\ZZ$ the subspace
 $\U^{(r)}$ is the linear span of the elements of $\U(R)$ whose degrees do not exceed $r$.
Let
\ben
\gr\U(R)=\bigoplus_{r\in\ZZ} \U^{(r)}/\U^{(r-1)}
\een
be
the associated graded algebra. It inherits a $\CC[h]$-module structure from $\U(R)$.
Let $\overline{l}_{ij}^{\ts(-r)}$ denote
the image of   $l_{ij}^{\tss(-r)}$ in the $r$-th component of $\gr\U(R)$. We will also write
\ben
\overline{L}^{\ts +} (u)=\sum_{i,j=1}^n e_{ij}\otimes
\overline{l}_{ij}^{\ts +}(u)\ts\in\ts \End\CC^n \otimes \gr\U(R)[[u]],
\een
where
\ben
\overline{l}_{ij}^{\ts +}(u)=\delta_{ij}-h\tss\sum_{r=1}^\infty
\overline{l}_{ij}^{\ts(-r)}u^{r-1}.
\een

\bpr\label{prop:isomgraded}
We have an isomorphism
\ben
\gr \U(R)\otimes_{\CC[h]}\CC[[h]]\cong \Y^+(\gl_n)
\een
defined on the generators by
\beql{iso}
\overline{l}_{ij}^{\ts(-r)}\mapsto {t}_{ij}^{(-r)}
\eeq
for all $1\leqslant i,j\leqslant n$ and $r=1,2,\dots$.
\epr

\bpf
Observe that both sides of
\eqref{RTT} are elements of finite degrees of
the algebra
\ben
\End\CC^n \ot\End\CC^n\ot\U(R)((u))((v))[[h]],
\een
with respect to
the degree function defined by \eqref{filtr} and \eqref{filtr2}
together with $\deg v=-1$. Moreover,
by taking
the highest degree components (which are of the zero degree)
we get the defining relations \eqref{RTTy} for the dual Yangian $\Y^+(\gl_n)$:
\ben
R^{\ts\textrm{rat}}(u-v)\tss\overline{L}_1^{\ts +}(u)\tss\overline{L}_2^{\ts +}(v)
=\overline{L}_2^{\ts +}(v)\tss\overline{L}_1^{\ts +}(u)\tss R^{\ts\textrm{rat}}(u-v).
\een
Note that the coefficients of any monomial $u^a v^b$ on both sides
of this relation coincide
with the highest degree components of the respective coefficients of the monomial $u^a v^b$
on both sides of \eqref{RTT}.
This implies that the mapping
${t}_{ij}^{(-r)}\mapsto\overline{l}_{ij}^{\ts(-r)}$ defines a homomorphism from
the Yangian to the extended graded algebra.
This homomorphism is clearly surjective,
while its injectivity follows from well-known versions of the
Poincar\'{e}--Birkhoff--Witt theorem
for the algebras $\Y^+(\gl_n)$
and $\U(R)$; see~\cite[Section~3.4]{ek:ql3}.
\epf

\section{Quantum affine vertex algebra}
\label{sec:QVOA}

We follow \cite{ek:ql5} to introduce a quantum vertex algebra structure on the
$h$-adic completion $\wt\U(R)$ of the quantized universal enveloping algebra $\U(R)$.
The corresponding structure associated with the rational $R$-matrix was studied in
\cite{jkmy:cq} where a description of the center of the quantum vertex algebra was given.
Our goal is to obtain an analogous description of the center of $\wt\U(R)$
at the critical level.

\subsection{The center of a quantum vertex algebra}
\label{subsec:QVA}

We shall say that the $\CC[[h]]$-module $V$ is topologically
free if $V=V_0[[h]]$ for some complex vector space $V_0$. Denote by $V_h ((z))$
the space of all Laurent series
\ben
v(z)=\sum_{r\in\mathbb{Z}} v_r z^{-r-1} \in V[[z^{\pm 1}]]
\een
satisfying $v_r\to 0$ as $r\to\infty$, in the $h$-adic topology.

\bde\label{def:qvoa}
Let $V=V_0 [[h]]$ be a topologically free $\mathbb{C}[[h]]$-module.
A {\em quantum vertex algebra} $V$ over $\mathbb{C}[[h]]$ is the following data.
\begin{enumerate}[(a)]
\item\label{101} A $\mathbb{C}[[h]]$-module map (the {\em vertex operators})
\beql{truncation}
Y\,\colon\, V\otimes V \,\to\, V_h ((z)),\qquad v\ot w\mapsto Y(z)\tss (v\ot w).
\eeq
Setting $Y(v,z)\tss w =\, Y(z)\tss (v\ot w)$
defines the map $Y(v,z):V\to V_h ((z))$
which satisfies the {\em weak associativity property}:
for any $u,v,w\in V$ and $n\in\mathbb{Z}_{\geqslant 0}$
there exists $\ell\in\mathbb{Z}_{\geqslant 0}$
such that
\ben
(z_0 +z_2)^\ell\ts Y(v,z_0 +z_2)Y(w,z_2)\tss u - (z_0 +z_2)^\ell\ts Y\big(Y(v,z_0)w,z_2\big)\tss u
\in h^n V[[z_0^{\pm 1},z_2^{\pm 1}]].
\een
\item\label{102} A vector $\vac\in V$ (the {\em vacuum vector}) which satisfies
$Y(\vac ,z)v=v$ for all $v\in V$,
and for any $v\in V$ the series $Y(v,z)\tss\vac$ is a Taylor series in $z$ with
the property
\beql{v2}
Y(v,z)\tss\vac\big|^{}_{z=0} =v.
\eeq
\item\label{103} A $\mathbb{C}[[h]]$-module map $D\colon V\to V$
(the {\em translation operator}) which satisfies
\ben
D\tss\vac =0\Fand
\frac{d}{dz}Y(v,z)=[D,Y(v,z)]\quad\text{for all }v\in V.
\een
\item\label{104} A $\mathbb{C}[[h]]$-module map
$\mathcal{S}\colon V\otimes V\to V\otimes V\otimes\mathbb{C}((z))$ which satisfies
\begin{align}
&\mathcal{S}(z)(v\otimes w)-v\otimes w\ot 1 \,\in\, h \,V\otimes V\otimes
\mathbb{C}((z))\quad \text{for } v,w\in V,\non\\
&[D\otimes 1, \mathcal{S}(z)]=-\frac{d}{dz}\mathcal{S}(z),\non\\
\intertext{the {\em Yang--Baxter equation}}
&\mathcal{S}_{12}(z_1)\tss\mathcal{S}_{13}(z_1+z_2)\tss\mathcal{S}_{23}(z_2)
=\mathcal{S}_{23}(z_2)\tss\mathcal{S}_{13}(z_1+z_2)\tss\mathcal{S}_{12}(z_1),\non
\end{align}
the {\em unitarity condition} $\mathcal{S}_{21}(z)=\mathcal{S}^{-1}(-z)$,
and the $\mathcal{S}$-{\em locality}:
for any $v,w\in V$ and $n\in\mathbb{Z}_{\geqslant 0}$ there exists
$\ell\in\mathbb{Z}_{\geqslant 0}$ such that for any $u\in V$
\begin{align}
&(z_1-z_2)^{\ell}\ts Y(z_1)\big(1\otimes Y(z_2)\big)\big(\mathcal{S}(z_1 -z_2)(v\otimes w)\otimes u\big)
\nonumber\\[0.4em]
&\qquad-(z_1-z_2)^{\ell}\ts Y(z_2)\big(1\otimes Y(z_1)\big)(w\otimes v\otimes u)
\in h^n V[[z_1^{\pm 1},z_2^{\pm 1}]].
\non
\end{align}
\end{enumerate}
\vspace{-2.0em}
\qed
\ede

The
tensor products in Definition~\ref{def:qvoa} are understood as
$h$-adically completed. In particular, $V\otimes V$ denotes the space
$(V_0 \otimes V_0)[[h]]$
and $V\otimes V\otimes\mathbb{C}((z))$ denotes the space
$\big(V_0\otimes V_0\otimes\mathbb{C}((z))\big)[[h]]$.

Let $V$ be a quantum vertex algebra. As in \cite{jkmy:cq}, we
define the {\em center} of $V$ as the $\CC[[h]]$-submodule
\ben
\z(V)=\big\{v\in V \ |\  w_rv=0\text{ for all }w\in V\text{ and all }r\geqslant 0 \big\}.
\een
It was proved in \cite{jkmy:cq} that the center of a quantum vertex
algebra is a unital associative algebra with the product $\z(V)\ot\z(V)\to\z(V)$
given by
\ben
v\cdot w=v_{-1}w\quad\text{for all }v,w\in V.
\een
The algebra $\z(V)$ need not be commutative; see \cite[Proposition~4.3]{jkmy:cq}.
Instead, it possesses the
$\Sc$-{\em commutativity} property as demonstrated in \cite[Proposition~3.7]{jkmy:cq}.
The next proposition shows that this property (as given in \eqref{scomm} below)
is characteristic for elements of the center.

\bpr\label{prop:scom}
Let $V$ be a quantum vertex algebra. Vector $v\in V$ belongs to $\z(V)$ if and only if
\beql{scomm}
Y(z_1)\left(1\ot Y(z_2)\right)\left(\Sc(z_1 -z_2)(v\ot w)\ot u  \right)
=Y(w,z_2)Y(v,z_1)\tss u
\eeq
for all $w\in V$ and $u\in\z(V)$.
\epr

\bpf
Let $v\in V$ satisfy \eqref{scomm} for all $w\in V$
and $u\in\z(V)$. Recall the $\Sc$-Jacobi identity:
\begin{align}
&z_{0}^{-1}\delta\left(\frac{z_2 -z_1}{z_0}\right)Y(w,z_2)Y(v,z_1)\tss u\nonumber\\
&\quad-z_{0}^{-1}\delta\left(\frac{z_1 -z_2}{-z_0}\right)
Y(z_1)\big(1\otimes Y(z_2)\big)\big(\mathcal{S}(-z_0)(v\otimes w)\otimes u\big)\nonumber\\
&\quad\quad=z_1^{-1} \delta\left(\frac{z_2 -z_0}{z_1}\right)
Y\big(Y(w,z_0)v,z_1\big)\tss u,\non
\end{align}
which holds in any quantum vertex algebra; see \cite{l:hq}.
Taking the residue   $\Res_{z_0}$ on both sides,
we obtain the $\Sc$-commutator formula
\begin{align*}
Y(w,z_2)Y(v,z_1)\tss u
-Y(z_1)\left(1\ot Y(z_2)\right)\left(\Sc(z_1 -z_2)(v\ot w)\ot u  \right)&\\
\quad =
\Res_{z_0} z_1^{-1} \delta\left(\frac{z_2 -z_0}{z_1}\right)
Y\big(Y(w,z_0)v,z_1\big)\tss u.&
\end{align*}
The left hand side is equal to zero by \eqref{scomm}, so that
\ben
\Res_{z_0} z_1^{-1} \delta\left(\frac{z_2 -z_0}{z_1}\right)
Y\big(Y(w,z_0)v,z_1\big)\tss u=0.
\een
This implies
\ben
\sum_{r\geqslant 0} (-1)^r \binom{b+r}{r} (w_r v)_{-a-b-r-2}
u=0\qquad\text{for all}\quad a,b\in\mathbb{Z}.
\een
Let $m> 0$. By \eqref{truncation} there exists $r_0>0$ such
that $w_r v \in h^m V$ for all $r>r_0$. Therefore,
\ben
\sum_{r= 0}^{r_0} (-1)^r \binom{b+r}{r} (w_r v)_{-a-b-r-2} u=0\mod
h^mV \quad\text{for all}\quad a,b\in\mathbb{Z}.
\een
By evaluating $(a,b)=(-i-1+c,i)$ for all $i=0,\dots,r_0$ with a fixed integer $c$,
we obtain a system of
$r_0 +1$ homogeneous linear equations in the variables $(w_r v)_{-r-c-1}u$ with $r=0,\dots,r_0$.
It is easily verified that its matrix
is invertible, so there is a unique solution,
\ben
(w_r v)_{-r-c-1}u=0\mod h^m V\qquad\text{for all}\quad r=0,\dots,r_0.
\een
By taking here $c=-r$, $u=\vac$ and using \eqref{v2} we get $w_r v \in h^m V$
for all $r=0,\dots,r_0$. We may conclude that
$w_r v \in h^m V$ for all $r,m\geqslant 0$.  Since $V$ is a topologically
free $\CC[[h]]$-module, $V$ is separated, so that $\cap_{m\geqslant 1}h^m V=0$.
This implies that $w_r v=0$ for all $r\geqslant 0$  which means that
the vector $v$ belongs to the center of $V$.

The ``only if{\ts}" part holds due to \cite[Proposition 3.7]{jkmy:cq}.
\epf

\subsection{The center of the quantum affine vertex algebra}
\label{subsec:center}

Here we recall some results of \cite{ek:ql5} describing
the quantum vertex algebra structure on $\U(R)$.
Let $\wt{\U}(R)=\U(R)[[h]]$ be the $h$-adic completion of $\U(R)$.
The following property will play a central role; see
\cite[Lemma 2.1]{ek:ql5}. For any nonnegative integer $m$,
consider the tensor product space
$
(\End\mathbb{C}^{n})^{\otimes {(m+1)}} \otimes \wt{\U}(R)
$
with the copies of the endomorphism algebra labelled by $0,1,\dots,m$.
Let $v_1,\dots,v_m$ be variables.

\ble\label{lem:EK}
For any $c\in\CC$ there exists a unique series
$L(u)\in \End\mathbb{C}^{n}\ot\big(\End \wt{\U}(R)\big)_h((u))$
such that for all $m\geqslant 0$ we have
\begin{multline}
L_{0}(u) L_{1}^+ (v_1)\dots L_{m}^+ (v_m)1 =
R_{0\tss 1}(u-v_1+hc/2)^{-1}\dots R_{0\tss m}(u-v_m+hc/2)^{-1}\\[0.5em]
{}\times\tss L_{1}^+ (v_1)\dots L_{m}^+ (v_m)
R_{0\tss m}(u-v_m-hc/2)\dots R_{0\tss 1}(u-v_1-hc/2)1.
\non
\end{multline}
\ele

Fix an arbitrary  complex number $c$. As
the action of the operator $L(u)$ on $\wt{\U}(R)$
depends on the choice of $c$, we will indicate this dependence by denoting
the completed quantized universal enveloping algebra $\wt{\U}(R)$ by $\wt{\U}_c(R)$.
The complex number
$c$ will be called the {\em level} of $\wt{\U}_c(R)$.
This terminology is motivated by the fact
that the classical limit
of the quantum vertex algebra $\wt{\U}_c(R)$ coincides with the affine vertex
algebra for $\gl_n$ at the level $c$; see~\cite{ek:ql5}.
The following relations hold for operators on $\End\mathbb{C}^{n} \otimes
\End\mathbb{C}^{n} \otimes \wt{\U}_c(R)$:
\begin{align}
R(u-v)L_1(u)L_2(v)&=L_2(v)L_1(u)R(u-v),\label{RTT2}\\[0.5em]
R(u-v+h\tss c/2)L_1(u)L_2^+(v)&=L^+_2(v)L_1(u)R(u-v-h\tss c/2).\label{RTT3}
\end{align}
Given a variable $z$ and a
family of variables
$u=(u_1,\dots,u_m)$, set
\ben
L^{}_{[m]}(u|z)=L_{1}(z+u_1)\dots L_{m}(z+u_m),\qquad
L_{[m]}^{+}(u|z)=L_{1}^{+}(z+u_1)\dots L_{m}^{+}(z+u_m).
\een
The respective components of the matrices $L^{+}(u)$ and $L(u)$
are understood as operators on $\wt{\U}_c(R)$.
The series $L_{i}(z+u_i)$
should be expanded in  nonnegative powers
of $u_i$.

By the results of
Etingof and Kazhdan~\cite{ek:ql5},
for any $c\in\CC$ there exists a unique well-defined structure of quantum vertex algebra
on the quantized universal enveloping algebra $\wt{\U}_{c}(R)$. In particular,
the vacuum vector is
$
\vac=1\in \wt{\U}_{c}(R),
$
the vertex operators are defined by
\ben
Y\big(L_{[m]}^{+}(u|0)\vac,z\big)=L_{[m]}^{+}(u|z)\ts L^{}_{[m]} (u|z+h\tss c/2)^{-1}
\een
and the translation operator $D$ is given by
\ben
e^{zD}\ts L^+ (u_1)\dots L^+ (u_m)\vac = L^+ (z+u_1)\dots L^+ (z+u_m)\vac.
\een
We will not reproduce
the definition of the map
$\mathcal{S}$ as it requires some additional notation and it will not be used below.
The map is defined in the same way as for the rational $R$-matrix;
see also \cite[Section~4.2]{jkmy:cq}.

Consider the $h$-permutation operator $P^h\in\End\CC^n\ot\End\CC^n[[h]]$ defined by
\beql{ph}
P^h = \sum_{i}e_{ii}\ot e_{ii} +e^{h/2}\sum_{i>j} e_{ij}\ot e_{ji}
+ e^{-h/2}\sum_{i<j}e_{ij}\ot e_{ji}.
\eeq
The symmetric group $\Sym_k$ acts on the space $(\CC^n)^{\ot k}$ by
$s_a\mapsto P_{s_a}^h=P_{a\ts a+1}^h$ for $a=1,\dots,k-1$, where $s_a$
denotes the transposition $(a,a+1)$. For a reduced decomposition
$\sigma=s_{a_1}\dots s_{a_l}$ of an element $\sigma\in\Sym_k$
set $P_{\sigma}^h =P^h_{s_{a_1}}\dots P^h_{s_{a_l}}$. Denote
by $A^{(k)}$ the image of the normalized symmetrizer under this action, so that
\ben
A^{(k)}=\frac{1}{k!}\sum_{\sigma\in\Sym_k}\sgn\sigma\cdot P^{h}_\sigma.
\een

Using the two-parameter $R$-matrix \eqref{R2p},
for arbitrary variables $u_1,\dots,u_k$ set
\ben
\overline{R}(u_1,\dots,u_k)=\prod_{1\leqslant a<b\leqslant k}
\overline{R}_{ab}(u_a,u_b),
\een
where the product is taken in the lexicographical order on the pairs $(a,b)$.
Note that \eqref{RTT} implies the relation
\beql{comm1}
\overline{R}(u_1,\dots,u_k)L_{1}^{+}(u_1)\dots  L_{k}^{+}(u_k)
=L_{k}^{+}(u_k)\dots L_{1}^{+}(u_1) \overline{R}(u_1,\dots,u_k).
\eeq
We will need the following well-known case of the fusion
procedure for the $R$-matrix \eqref{R2p} going back to \cite{c:ni}.

\ble\label{fusionl}
Set $u_a = u-(a-1)h$ for $a=1,\dots,k$. We have
\ben
\overline{R}(u_1,\dots,u_k)=k! \ts e^{uk(k-1)/2}
\prod_{0\leqslant a<b\leqslant k-1}(e^{-ah}-e^{-bh})A^{(k)}.
\vspace{-1.5em}
\een
\qed
\ele

Combining relation \eqref{comm1} and Lemma~\ref{fusionl} we get for $u_a = u-(a-1)h$
\beql{AT}
A^{(k)}L_{1}^{+}(u_1)\dots L_{k}^{+}(u_k)=L_{k}^{+}(u_k)\dots L_{1}^{+}(u_1) A^{(k)}.
\eeq
Also, since
\ben
\overline{R}(u,v)D_1 D_2 = D_2 D_1 \overline{R}(u,v),
\een
we have
\beql{AD}
A^{(k)}D_1\dots D_k =D_k\dots D_1 A^{(k)}
\eeq
for the diagonal matrix \eqref{diagonal}.

Now consider the $h$-adically completed
quantum vertex algebra $\wt{\U}_{c}(R)$ at the critical level $c=-n$
and introduce its elements as the coefficients of the power series in $u$ given by
\beql{phik}
\phi_k(u)=\tr_{1,\dots,k} A^{(k)} L_{1}^+ (u) \dots L_{k}^{+}(u-(k-1)h)D_1 \dots D_k
\eeq
for $k=1,\dots,n$.

\bpr\label{prop:lkcen}
The coefficients of all series $\phi_k(u)$ belong to
$\z(\wt{\U}_{-n}(R))$.
\epr

\bpf
The argument is essentially a version of the proofs of \cite[Theorem~3.2]{fjmr:hs}
and \cite[Theorem~2.4]{jkmy:cq} which we adjust for the current context.
It is sufficient to show that
\ben
L_0(v)\tss\phi_k(u)=\phi_k(u)
\een
for all $k=1,\dots,n$.
Relation
\eqref{RTT3} implies
\begin{multline}
L_0(v) L_{1}^+ (u_1) \dots L_{k}^{+}(u_k)
=R_{01}(v-u_1-hn/2)^{-1}\dots R_{0k}(v-u_k-hn/2)^{-1}\\[0.5em]
{}\times\ts L_{1}^+ (u_1) \dots L_{k}^{+}(u_k)L_0(v)
R_{0k}(v-u_k+hn/2)\dots R_{01}(v-u_1+hn/2).\label{c1}
\end{multline}
Since $L_0(v)$ commutes with $A^{(k)}$ and $D_a$ for $a=1,\dots,k$,
using \eqref{c1} and $L_0(v)\vac = \vac$ we get
\begin{multline}
L_0(v)\tss \phi_k(u)=
\tr_{1,\dots,k}\ts A^{(k)} R_{01}(v-u_1-hn/2)^{-1}\dots R_{0k}(v-u_k-hn/2)^{-1}\\[0.5em]
{}\times\ts L_{1}^+ (u_1) \dots L_{k}^{+}(u_k)
 \ts R_{0k}(v-u_k+hn/2)\dots R_{01}(v-u_1+hn/2)D_1 \dots D_k.\label{c2}
\end{multline}
Hence, we only have to prove that  the right hand side in \eqref{c2}  equals $\phi_k(u)$.
Set
\begin{align}
X={}&R_{01}(v-u_1-hn/2)^{-1}\dots R_{0k}(v-u_k-hn/2)^{-1}\non\\[0.4em]
Y={}&L_{1}^+ (u_1) \dots L_{k}^{+}(u_k)   R_{0k}(v-u_k+hn/2)
\dots R_{01}(v-u_1+hn/2)D_1 \dots D_k\non
\end{align}
and denote by $X^{\textrm{op}}$ and $Y^{\textrm{op}}$ the respective expressions
obtained from these products by reversing the order of some factors:
\begin{align*}
X^{\textrm{op}}={}&R_{0k}(v-u_k-hn/2)^{-1}\dots R_{01}(v-u_1-hn/2)^{-1} \\[0.4em]
Y^{\textrm{op}}={}&L_{k}^{+}(u_k)  \dots L_{1}^+ (u_1) R_{01}(v-u_1+hn/2)
\dots  R_{0k}(v-u_k+hn/2)D_k \dots D_1.
\end{align*}
The Yang-Baxter equation for the $R$-matrix \eqref{R2p} and Lemma \ref{fusionl} imply
\begin{align}
&A^{(k)}R_{01}(v-u_1-hn/2)^{-1}\dots R_{0k}(v-u_k-hn/2)^{-1}\nonumber\\[0.4em]
&\qquad= R_{0k}(v-u_k-hn/2)^{-1}\dots R_{01}(v-u_1-hn/2)^{-1} A^{(k)},
\non
\end{align}
that is,
$
A^{(k)}X=X^{\textrm{op}}A^{(k)}$,
and also
\begin{align}
&A^{(k)}R_{0k}(v-u_k+hn/2)\dots R_{01}(v-u_1+hn/2)\nonumber\\[0.4em]
&\qquad= R_{01}(v-u_1+hn/2)\dots R_{0k}(v-u_k+hn/2) A^{(k)}.
\label{c7}
\end{align}
Combining \eqref{AT}, \eqref{AD} and \eqref{c7} we get
$
A^{(k)}Y=Y^{\textrm{op}}A^{(k)}$.
Since $A^{(k)}$ is an idempotent, we proceed as follows:
\begin{align*}
\tr_{1,\dots,k} A^{(k)}XY
&=\tr_{1,\dots,k} X^{\textrm{op}}A^{(k)}Y
=\tr_{1,\dots,k} X^{\textrm{op}}\left(A^{(k)}\right)^2 Y\\[0.4em]
&=\tr_{1,\dots,k} X^{\textrm{op}}A^{(k)}A^{(k)} Y
=\tr_{1,\dots,k} A^{(k)}X Y^{\textrm{op}} A^{(k)}.
\end{align*}
By the cyclic property of trace, this equals
\ben
\tr_{1,\dots,k} A^{(k)}X Y^{\textrm{op}} A^{(k)}
=\tr_{1,\dots,k} X Y^{\textrm{op}} \left(A^{(k)}\right)^2
=\tr_{1,\dots,k} X Y^{\textrm{op}} A^{(k)}
=\tr_{1,\dots,k} X A^{(k)}  Y,
\een
and so
$
\tr_{1,\dots,k} A^{(k)}XY=\tr_{1,\dots,k} X A^{(k)}  Y.
$
As a final step, we use the property
\ben
\tr_{1,\dots,k} X A^{(k)}  Y=\tr_{1,\dots,k} X^{t_1 \dots t_k} \left(A^{(k)}  Y\right)^{t_1 \dots t_k}.
\een
Write
\begin{align*}
&X^{t_1 \dots t_k} \left(A^{(k)}Y\right)^{t_1 \dots t_k}
=\left(R_{01}(v-u_1-hn/2)^{-1}\right)^{t_1}\dots \left(R_{0k}(v-u_k-hn/2)^{-1}\right)^{t_k}\\[0.5em]
&\quad\times\tss D_1\dots D_k R_{0k}(v-u_k+hn/2)^{t_k}\dots
R_{01}(v-u_1+hn/2)^{t_1}\left(A^{(k)}L_{1}^+ (u_1) \dots L_{k}^{+}(u_k)\right)^{t_1 \dots t_k}
\end{align*}
and apply the first crossing symmetry property in \eqref{csym} to get
\begin{multline}
X^{t_1 \dots t_k} \left(A^{(k)}Y\right)^{t_1 \dots t_k}
= D_1\dots D_k \left(A^{(k)}L_{1}^+ (u_1) \dots L_{k}^{+}(u_k)\right)^{t_1 \dots t_k}\\[0.4em]
{}= \left(A^{(k)}L_{1}^+ (u_1) \dots L_{k}^{+}(u_k)D_1\dots D_k\right)^{t_1 \dots t_k}.\non
\end{multline}
This implies
\begin{align*}
\tr_{1,\dots,k} X A^{(k)}  Y&=\tr_{1,\dots,k} \left(A^{(k)}L_{1}^+ (u_1)
\dots L_{k}^{+}(u_k)D_1\dots D_k\right)^{t_1 \dots t_k}\\[0.4em]
 &=\tr_{1,\dots,k} A^{(k)}L_{1}^+ (u_1) \dots L_{k}^{+}(u_k)D_1\dots D_k=\phi_k(u),
 \end{align*}
thus completing the proof.
\epf

Consider the $h$-permutation operator $P^h_{(k,k-1,\dots,1)}$ which is
associated with the $k$-cycle $(k,k-1,\dots,1)=s_{k-1}\dots s_1$ so that
\ben
P^h_{(k,k-1,\dots,1)}=P^h_{k-1\ts k}\dots P^h_{1\tss 2}.
\een
Introduce another family of elements of $\wt{\U}_{-n}(R)$ as the coefficients of the power series in $u$
defined by
\beql{tracepow}
\theta_k(u)=\tr_{1,\dots,k} P^h_{(k,k-1,\dots,1)} L_{1}^+ (u) \dots
L_{k}^{+}(u-(k-1)h)\tss D_1 \dots D_k
\eeq
for all $k\geqslant 1$ and $\theta_0(u)=n$.

\bco\label{cor:trpow}
The coefficients of all series $\theta_k(u)$ belong to
$\z(\wt{\U}_{-n}(R))$.
\eco

\bpf
Relation \eqref{RTT} implies that
$M=L^+(u)\tss D\tss e^{-h\di_u}$
is a $q$-{\em Manin matrix} with entries in
$\wt{\U}_{-n}(R)[[u,\partial_u]]$ as follows from \cite[Lemma~5.1]{cfrs:ap},
with the notation \eqref{xq}.
Hence, applying the Newton identity for $q$-Manin
matrices \cite[Theorem~5.7]{cfrs:ap}, we can express the coefficients of
each series $\theta_k(u)$ as polynomials in the coefficients of the series
$\phi_k(u)$. So the corollary is a consequence of
Proposition~\ref{prop:lkcen}.
\epf

Let $\wt{\U}(R)^{\text{\rm ext}}$ denote the extension
of $\wt{\U}(R)$ to an algebra over the field $\CC((h))$. For all $m\geqslant 0$
introduce elements $\Theta_{m}^{(r)}$ of this algebra as the coefficients of the series
\ben
\Theta_m (u)=\sum_{r=0}^{\infty}\Theta_m^{(r)}u^r \in \wt{\U}(R)^{\text{\rm ext}}[[u]],
\een
where we use the series \eqref{tracepow} and set
\beql{Theta}
\Theta_m (u)=h^{-m}\sum_{k=0}^m (-1)^k\binom{m}{k}\tss
\theta_k(u).
\eeq

\bpr\label{prop:commtheta}
All the elements $\Theta_{m}^{(r)}$ with $r,m\geqslant 0$ pairwise commute.
\epr

\bpf
This follows from the corresponding well-known property of the coefficients of
the series \eqref{phik}; see e.g.~\cite[Proposition~6.5]{cfrs:ap} for a proof, which is
quite similar to the rational case; cf.~\cite[Section~1.14]{m:yc}. The property extends to
the coefficients of the series \eqref{tracepow} due to the Newton identity;
see \cite[Theorem~6.6]{cfrs:ap}.
\epf

\bth\label{thm:algindep}
All coefficients of the series $\Theta_{m}(u)$ belong to the
$\CC[[h]]$-module $\z(\wt{\U}_{-n}(R))$.
Moreover, the family  $\Theta_{m}^{(r)}$ with $m=1,\dots,n$
and $r=0,1,\dots$ is algebraically independent.
\eth

\bpf
We will need the usual permutation $P$ given in \eqref{permut}
along with the $h$-permutation operator $P^h$ defined in
\eqref{ph}. Let $M=L^+(u)\tss D\tss e^{-h\di_u}$ as before, and
for each $m\geqslant 1$ consider the expression
\begin{multline}\label{temp_a}
\Mc_m=h^{-m}\big(1-(M_m)^{\rightarrow}  \big)
\Big(P^{}_{m-1\,m}-P^h_{m-1\,m}(M_{m-1})^{\rightarrow}\Big)\\[0.4em]
{}\times\dots\times\Big(P^{}_{23}-P^h_{23}(M_2)^{\rightarrow} \Big)
\Big(P^{}_{12}-P^h_{12}M^{}_1  \Big),
\end{multline}
where the arrow in the superscript indicates that the corresponding factor appears on the right:
\beql{apptx}
\Big(P^{}_{a\ts a+1}-P^h_{a\ts a+1}\big(M_a \big)^{\rightarrow} \Big)\ts X:
=P^{}_{a\ts a+1}\ts X-P^h_{a\ts a+1}\ts X\ts M_a.
\eeq
We verify first that the expression \eqref{temp_a}, as a Laurent series in $h$,
does not contain negative powers of $h$.
Indeed, write
\ben
P^{}_{a\,a+1}-P^h_{a\,a+1}\big(M_a\big)^{\rightarrow}
=P^{}_{a\,a+1}\Big(1-P^{}_{a\,a+1}P^h_{a\,a+1}
\big(L_{a}^+(u)D^{}_{a} e^{-h\partial_u}\big)^{\rightarrow}  \Big)
\een
and observe that
\ben
P^{}_{a\,a+1}P^h_{a\,a+1}\big|_{h=0}=1.
\een
Therefore, the expression in \eqref{apptx} vanishes at $h=0$ for any element $X$
of the $\CC[[h]]$-module $\left(\End \CC^n\right)^{\ot m}\otimes \wt{\U}_{-n}(R)$.
Hence each of the $m$ factors in \eqref{temp_a} is divisible by $h$.

As a next step, expand the product in \eqref{temp_a} to get the expression
\ben
\Mc_m=h^{-m}\ts\sum_{k=0}^m \sum_{1\leqslant a_1<\dots< a_k\leqslant m }\ts(-1)^k\ts
\Pi_{a_1,\dots,a_k}\ts M_{a_1}\dots M_{a_k},
\een
where
\ben
\Pi_{a_1,\dots,a_k}=P^{}_{(m,m-1,\dots,a_k+1)}P^h_{a_k\ts a_k+1}
P^{}_{(a_k,\dots,a_{k-1}+1)}P^h_{a_{k-1}\ts a_{k-1}+1}\dots
P^{}_{(a_2,\dots,a_1+1)}P^h_{a_1\ts a_1+1}P^{}_{(a_1,\dots,1)}.
\een
Now consider the trace $\tr_{1,\dots,m}\tss\Mc_m$. Let us verify that for the partial trace
we have
\beql{trpi}
\tr^{}_{\{1,\dots,m\}\setminus\{a_1,\dots, a_k\}}\ts \Pi_{a_1,\dots,a_k}
=P^h_{a_{k-1}\ts a_{k}}P^h_{a_{k-2}\ts a_{k-1}}\dots P^h_{a_1\ts a_2}.
\eeq
For any permutation $\si\in\Sym_m$ we have
$P^{}_{\si}\tss P^h_{a\tss b}=P^h_{\si(a)\tss \si(b)}P^{}_{\si}$. Moreover,
$\tr_{a}P^{}_{a\tss b}=\tr_{a}P^h_{a\tss b}=1$ for any $a\ne b$. Therefore,
\ben
\tr_{a_{r}+1,\dots,a_{r+1}-1}\ts P^{}_{(a_{r+1},\dots,a_{r}+1)}\ts P^h_{a_{r}\ts a_{r}+1}
=\tr_{a_{r}+1,\dots,a_{r+1}-1}\ts P^h_{a_{r}\ts a_{r+1}}\ts P^{}_{(a_{r+1},\dots,a_{r}+1)}=
P^h_{a_{r}\ts a_{r+1}}
\een
for $r=1,\dots,k$ with $a_{k+1}:=m$. This implies \eqref{trpi}. Hence we obtain
\ben
\bal
\tr_{1,\dots,m}\tss\Mc_m&=
h^{-m}\ts\sum_{k=0}^m \sum_{1\leqslant a_1<\dots< a_k\leqslant m }
\ts(-1)^k\ts\tr^{}_{a_1,\dots,a_k}\ts
P^h_{a_{k-1}\ts a_{k}}P^h_{a_{k-2}\ts a_{k-1}}
\dots P^h_{a_1\ts a_2}\ts M_{a_1}\dots M_{a_k}\\[0.4em]
{}&=h^{-m}\ts\sum_{k=0}^m \ts(-1)^k\ts \binom{m}{k}\tss \tr_{1,\dots,k}\ts
P^h_{k-1\ts k}P^h_{k-2\ts k-1}\dots P^h_{1\ts 2}\ts M_1\dots M_k.
\eal
\een
Since $P^h_{k-1\ts k}P^h_{k-2\ts k-1}\dots P^h_{1\ts 2}=P^h_{(k,k-1,\dots,1)}$ and
\ben
M_1\dots M_k= L_{1}^+ (u) \dots
L_{k}^{+}(u-(k-1)h)\tss D_1 \dots D_k\tss e^{-kh\partial_u},
\een
using \eqref{tracepow} we can write
\ben
\tr_{1,\dots,m}\tss\Mc_m=h^{-m}\ts\sum_{k=0}^m \ts(-1)^k\ts \binom{m}{k}\tss
\theta_k(u)\tss e^{-kh\partial_u}.
\een
Regarding this expression as a power series in $\di_u$, we conclude that
its constant term coincides with $\Theta_m (u)$ as defined in \eqref{Theta}.
The first part of the theorem now follows from Corollary~\ref{cor:trpow}.

Note that the level $c$ is irrelevant for the second part of the theorem.
Extend the degree function \eqref{filtr2}
to the algebra $\wt{\U}(R)$
by allowing it to take the infinite value. Elements of finite degree
will then form a subalgebra which we denote by
$\wt{\U}(R)_{\text{fin}}$. Similarly, by using \eqref{filtr} introduce
the subalgebra
$\wt{\U}(R)[[u]]_{\text{fin}}$ of $\wt{\U}(R)[[u]]$
formed by elements of finite degree.
Observe that
the series $\Theta_{m}(u)$ belongs to $\wt{\U}(R)[[u]]_{\text{fin}}$.
Take the highest degree component of this series (that is, the associated element
of the graded algebra $\gr\wt{\U}(R)[[u]]_{\text{fin}}$) and
identify it with an element of the algebra $\Y^+(\gl_n)[[h,u]]$ via
an extension of the map \eqref{iso}.
This element
coincides with the series
\beql{newlabel2}
\overline{\Theta}_m (u)=h^{-m}\sum_{k=0}^m (-1)^k\binom{m}{k}
\tr\ts
T^+(u)\cdots T^+ (u-(k-1)h).
\eeq
Write
\ben
\overline{\Theta}_{m}(u)=\sum_{r=0}^\infty \overline{\Theta}_{m}^{(r)}u^r.
\een
By \cite[Proposition 4.7]{jkmy:cq}, the
coefficients $\overline{\Theta}_{m}^{\ts(r)}$ with $m=1,\dots,n$ and $r=0,1,\dots$
are algebraically independent elements of the $h$-adically
completed dual Yangian $\Y^+(\gl_n)$.
Note that  the coefficient of any power $u^a $ on the right hand side of \eqref{newlabel2}
coincides with the highest
degree component of the coefficient of  $u^a $ on the right hand side of \eqref{Theta}.
Hence, the corresponding coefficients  $\Theta_{m}^{(r)}$ are also
algebraically independent.
\epf

\bth\label{thm:topolgen}
The center at the critical level $\z(\wt{\U}_{-n}(R))$ is a
commutative algebra. It is topologically generated by the family
$\Theta_{m}^{(r)}$ with $m=1,\dots,n$ and $r=0,1\dots$.
\eth

\bpf
By Proposition~\ref{prop:commtheta} and Theorem~\ref{thm:algindep},
the coefficients $\Theta_{m}^{(r)}$ generate a commutative subalgebra of $\z(\wt{\U}_{-n}(R))$.
Let $X$ be an arbitrary element of  $ \z(\wt{\U}_{-n}(R))$.
We will prove by induction that for each $k\geqslant 0$ there exists a polynomial
\ben
Q\in\CC[\Theta_{m}^{(r)} ] [h],\qquad m=1,\dots , n\quad\text{and}\quad r=0,1,\dots
\een
such that   $X-Q\in h^k\tss \wt{\U}_{-n}(R)$.
The induction base is clear.
Suppose that the property  holds for some $k\geqslant 0$ so that
\beql{temp_a4}
X-Q\ts=\ts h^k X_k +h^{k+1} X_{k+1}+\dots\quad \text{with }X_s\in V_0,
\eeq
where we write $\wt{\U}_{-n}(R)=V_0 [[h]]$ for some complex vector space $V_0$.
By Lemma \ref{lem:EK} the operator $L(u)$ can be written as
\ben
L(u)=\sum_{i,j=1}^ne_{ij}\ot l_{ij}(u)\qquad\text{with}
\quad l_{ij}(u)=\delta_{ij}+h\sum_{r\in\ZZ}
\tilde{l}_{ij}^{\ts(r)}u^r\in \big(\End \wt{\U}_{-n}(R)\big)_h((u)).
\een
Since $X-Q$ belongs to $\z(\wt{\U}_{-n}(R))$, we have
\ben
\tilde{l}_{ij}^{\ts(r)}(X-Q)=0\qquad\text{for all}\quad i,j=1,\dots, n \fand r\in\ZZ.
\een
Thus, \eqref{temp_a4} implies
\beql{temp_a3}
\tilde{l}_{ij}^{\ts(r)}X_k \equiv 0\mod h,\qquad\text{hence}\quad X_k\in \z(\wt{\U}_{-n}(R))\mod h.
\eeq

Consider the symbol (the highest degree component)
of $X_k$ in the graded algebra
$\gr \U_{-n}(R)$. Its image
under the isomorphism of Proposition \ref{prop:isomgraded}
is an element $\overline{X}_k\in\Y^{+}(\gl_n)$. We will also regard
it as an element of the quantum vertex algebra
$\Vc_{\text{cri}}=\Y^{+}(\gl_n)[[h]]$ at the critical
level $c=-n$, associated with the double Yangian; see \cite[Theorem 4.1]{jkmy:cq}
for a precise definition of the quantum vertex algebra structure.
As shown in \cite[Section~4.4]{jkmy:cq}, the center
of this quantum vertex algebra coincides with the subspace
of invariants
\beql{temp_a2}
\z(\Vc_{\text{cri}})=
\{U\in \Vc_{\text{cri}}\ |\ t_{ij}^{\ts(r)}\tss U=0\qquad
\text{for $r\geqslant 1$\ \  and all\ \  $i,j$}\}.
\eeq
Here the operators $t_{ij}^{\ts(r)}$ are found as the coefficients of
the series
\ben
t_{ij}(u)=\delta_{ij}+h\tss\sum_{r=1}^{\infty}
t_{ij}^{(r)}u^{-r}\in \big(\End \Vc_{\text{cri}}\big)_h((u)),
\een
which are the entries of the matrix operator $T(u)$
\ben
T(u)=\sum_{i,j=1}^n e_{ij}\ot t_{ij}(u)
\een
uniquely determined by the relations
\begin{multline}
T_{0}(u)\tss T_{1}^+ (v_1)\dots T_{m}^+ (v_m)1 =
R_{01}^{\ts\text{rat}}(u-v_1-hn/2)^{-1}\dots R_{0m}^{\ts\text{rat}}(u-v_m-hn/2)^{-1}\\[0.5em]
{}\times\tss T_{1}^+ (v_1)\dots T_{m}^+ (v_m)
R_{0m}^{\ts\text{rat}}(u-v_m+hn/2)\dots R_{01}^{\ts\text{rat}}(u-v_1+hn/2)1
\label{temp_c6}
\end{multline}
with the notation as in Lemma~\ref{lem:EK}.
Similar to the proof of Theorem~\ref{thm:algindep}, extend the degree function
defined in \eqref{filtr} and \eqref{filtr2} by setting $\deg v_i=-1$ for all $i$.
As pointed out in Section \ref{sec:QUE},
the rational $R$-matrix $R^{\ts\text{rat}}(u)$ coincides
with the highest degree component
of the trigonometric $R$-matrix $R(u)$.
Write the expression $\Lc=L_{0}(u)\tss L_{1}^+ (v_1)\dots L_{m}^+ (v_m)1$
as a series in the monomials $u^a\tss v^{b_1}_1\dots v^{b_m}_m$
with coefficients in $\wt{\U}_{-n}(R)$.
The symbol (the highest degree component)
of $\Lc$
is an element of the graded algebra
$\gr \wt{\U}_{-n}(R)((u))[[v_1,\dots,v_m,h]]$.
Also writing this symbol as a series in the monomials $u^a\tss v^{b_1}_1\dots v^{b_m}_m$,
observe that if $a\leqslant 0$ then the coefficient of such a monomial in the symbol
coincides with the symbol of the coefficient of
the same monomial in the expansion of $\Lc$. This implies that
the image of the component of the symbol of $\Lc$ corresponding
to nonpositive powers of $u$
under an extension of the map \eqref{iso}
coincides with the right-hand side in \eqref{temp_c6}.
Together with \eqref{temp_a3} this shows that
for $r\geqslant 1$ and all $i,j\in\{1,\dots,n\}$
\ben
t_{ij}^{(r)}\ts\overline{X}_k \equiv 0\mod h,\qquad\text{hence}\quad
\overline{X}_k\in \z(\Vc_{\text{cri}})\mod h.
\een
It was proved in \cite[Theorem 4.8]{jkmy:cq} that the elements
$\overline{\Theta}_m^{\ts(r)}$ topologically generate $\z(\Vc_{\text{cri}})$.
Therefore,
$\overline{X}_k\equiv \overline{S}\mod h$ for some polynomial
\ben
\overline{S}\in\CC[\overline{\Theta}_m^{\tss(r)} ],\qquad m=1,
\dots,n\quad\text{and}\quad r=0,1,\dots.
\een
Note that $\deg \overline{\Theta}_m^{\tss(r)}=\deg \overline{\Theta}_m^{\tss(r)}\big|_{h=0}$
and so
\beql{temp_1e}\deg \overline{S} = \deg
\overline{S}\big|_{h=0}=\deg\overline{X}_k.
\eeq
Replace the  variables $\overline{\Theta}_{m}^{\tss(r)}$ in
$\overline{S}$
with the respective elements $\Theta_{m}^{(r)}$ to get a polynomial
\ben
S\in\CC[\Theta_{m}^{(r)}].
\een
Let us consider the difference $X_k -S$ and take its symbol
(the highest degree component) in the graded algebra $\gr\wt{\U}(R)_{\text{fin}}$.
It
belongs to $h\ts\gr\wt{\U}(R)_{\text{fin}}$. Therefore, we can conclude that
\ben
X_k - S \in X^{(1)}_k +h\tss\wt{\U}_{-n}(R)
\een
for some $X^{(1)}_k\in V_0$ whose degree is lower than that of the symbol.
Hence, we can write \eqref{temp_a4} in the form
\beql{temp_a7}
X-Q-h^k S\ts=\ts h^k X^{(1)}_k +h^{k+1}
X^{(1)}_{k+1}+\dots\quad \text{with }X^{(1)}_s\in V_0.
\eeq
Recall that $\deg l_{ij}^{\tss(-r)}=r$, so the elements
of $V_0$ have nonnegative degrees; see \eqref{filtr2}.
Due to \eqref{temp_1e} we have
$\deg S = \deg X_k$, so that $\deg X_k>\deg X_k ^{(1)}\geqslant 0$.

Now we can repeat the same argument as above, but working with \eqref{temp_a7} instead of
\eqref{temp_a4} so that the role of $X_k$ is played by $X^{(1)}_k$.
An obvious induction on the degree of the element $\deg X_k$
allows us to conclude
that there exists a polynomial
$P\in\CC[\Theta_{m}^{(r)} ][h] $ satisfying
$
   X-P\in h^{k+1}\tss\wt{\U}_{-n}(R)
$
   thus completing the induction step and the proof.
\epf

Now consider the quantum vertex algebra $\wt{\U}_c (R)$ at a noncritical
level $c\neq -n$. Define the quantum determinant of the matrix $L^+ (u)$ by
\beql{qdet}
\qdet L^+ (u)=\sum_{\sigma\in\mathfrak{S}_n}
\left(-e^{-h/2} \right)^{l(\sigma)} l_{\sigma(1)1}^{+}(u)\dots l_{\sigma(n)n}^{+}(u-(n-1)h),
\eeq
where $l(\sigma)$ equals the number of inversions
in the sequence $(\sigma(1),\dots,\sigma(n))$.
Write
\ben
\qdet L^+ (u)=1-h\tss(d_0 +d_1 u +d_2 u^2+\dots).
\een

\bpr
The center   $\z(\wt{\U}_{c}(R))$ at a noncritical
level $c\neq -n$ is a commutative algebra. It is topologically generated by the family
$d_0,d_1,\dots$ of algebraically independent elements.
\epr

\bpf
The property that the coefficients $d_0,d_1,\dots$ are pairwise commuting elements of
$\z(\wt{\U}_{c}(R))$   can be verified by repeating the
arguments of the proof of \cite[Lemma~4.3]{fjmr:hs} with the use of Corollary~\ref{cor:gu}
(the assumption that the level is critical is unnecessary for the lemma to hold true;
see also \cite[Proposition~2.8]{jkmy:cq}).
In fact, these arguments demonstrate that the coefficients
belong to the center of the algebra $\wt{\U}_{c}(R)$ and so, in particular, they
commute pairwise.

Note that the series $\qdet L^+ (u)$ belongs to the algebra $ \wt{\U}(R)[[u]]_{\text{fin}}$
introduced in the proof of Theorem~\ref{thm:algindep}.
The symbol of $\qdet L^+ (u)$ (the highest degree component) belongs to the
graded algebra $\gr \wt{\U}(R)[[u]]_{\text{fin}}$. Identify the symbol
with an element of the $\Y^+ (\gl_n)[[u,h]]$ via an extension of the map \eqref{iso}.
This element
equals
\beql{qdetY}
\qdet T^+ (u)=\sum_{\sigma\in\mathfrak{S}_n} \sgn \sigma \cdot t_{\sigma(1)1}^{+}(u)\dots
t_{\sigma(n)n}^{+}(u-(n-1)h)
\in \Y^+ (\gl_n)[[u,h]].
\eeq
Now we will use some results concerning the quantum vertex algebra
$\Vc_{\text{c}}:=\Y^+ (\gl_n)[[h]]$ at the level $c\neq -n$;
see \cite[Theorem 4.1]{jkmy:cq}.
It was proved in \cite[Proposition 4.5]{jkmy:cq} that the center
$\z(\Vc_{\text{c}})$ is topologically generated by the algebraically
independent family of coefficients $\overline{d}_0,\overline{d}_1,\dots$ of the quantum determinant
\ben
\qdet T^+ (u)=1-h\tss(\overline{d}_0 +\overline{d}_1 u +\overline{d}_2 u^2+\dots).
\een
This implies that the coefficients $d_0,d_1,\dots$ are algebraically independent.
Finally, the property that the family $d_0,d_1,\dots$ topologically
generates the center $\z(\wt{\U}_{c}(R))$
is verified by the same argument as in the proof of Theorem~\ref{thm:topolgen}.
\epf

\appendix
\section{A direct proof of Proposition~\ref{prop:normal}}
\label{app:gu}

We start by describing solutions of
\eqref{g1}, regarding it as an equation for $g(u)\in\CC((u))[[h]]$
(we suppress the dependence on $n$ and $h$ in the notation
of these series).
Write
\beql{guexp}
g(u)=\sum_{l\geqslant 0}g_{l}(u) h^l,\qquad g_l (u)\in\CC((u)).
\eeq
Using the Taylor expansion formula
\ben
g(u+nh)=\sum_{l\geqslant 0}\Big(\sum_{k=0}^l \frac{n^k}{k!}g_{l-k}^{(k)}(u)\Big)h^l,
\een
from \eqref{g1} we get
\begin{align}
&\Bigg((1-e^u)-e^u (1-e^u)\sum_{l\geqslant 0 }\frac{n^l}{l!}h^l\Bigg)\cdot
\sum_{l\geqslant 0}\Big(\sum_{k=0}^l \frac{n^k}{k!}g_{l-k}^{(k)}(u)\Big)h^l\nonumber\\
&\qquad\qquad=\Bigg(1+e^u\sum_{l\geqslant 0}\frac{e^u n^l -(n-1)^l -1}{l!}h^l   \Bigg)
\cdot  \sum_{l\geqslant 0}g_{l}(u) h^l.\label{g2}
\end{align}
By equating the coefficients of the same powers of $h$ on  both sides of
\eqref{g2}, we get a system of differential equations
in $g_k(u)$ with $k\geqslant 0$. For the constant term the equation holds identically, while
considering the coefficients of $h$ in \eqref{g2}
we get  $g_{0}'(u)=0$ so that
$g_{0}(u)=c_0$ for $c_0\in\CC$.
Taking the coefficients of $h^l$ in \eqref{g2}, for an arbitrary $l\geqslant 1$ we get
\beql{g3}
g_{l-1}'(u)=\frac{e^u}{(1-e^u)^2 }\sum_{k=0}^{l-2}
\sum_{m=0}^{l-k-1}p_{k,m}(e^u)\tss g_{k}^{(m)}(u)-\sum_{k=0}^{l-2}
\frac{n^{l-k-1}}{(l-k)!}\tss g_{k}^{(l-k)}(u)
\eeq
for some polynomials  $p_{k,m}(z)\in\CC[z]$ of degree not exceeding $1$
such that
\ben
p_{k,m}(1)=0\quad\text{when }k+m=l-1,
\een
and $p_{k,0}(z)$ are constants. The right hand side of \eqref{g3} is understood as
being equal to zero for $l=1$.

Fix $r\geqslant 1$ and suppose that series $g_{0}(u),\dots,g_{r-1}(u)\in\CC((u))$
satisfy \eqref{g3} for $l=1,\dots,r$. We will prove by induction on $r$ that
a solution $g_r(u)\in \CC((u))$ of
\eqref{g3} with $l=r+1$
exists and, up to an additive constant, it has the form of a linear combination of the series
\ben
  \frac{ p(e^u)}{(1-e^{u})^t}
  \een
for some polynomials
$p(z)\in\CC[z]$ and positive integers $t$ such that
$t- \deg p\geqslant 1$ and $t\leqslant r$.
Indeed, by the induction hypothesis, all derivatives
$g^{(m)}_{k}(u)$
with $m\geqslant 1$ and $k=0,\dots, r-1$
can be expressed as linear combinations of series of the form
\ben
\frac{e^{u}\tss q(e^u)}{(1-e^{u})^{t+m}}
\een
for some polynomials
$q(z)\in\CC[z]$ and positive integers $t$ such that
$t+m- \deg q\geqslant 2$ and $t\leqslant k$.
This implies that the right hand side in \eqref{g3} for $l=r+1$ can
be written as a linear combination of series of the form
\ben
\frac{e^{u}\tss q(e^u)}{(1-e^{u})^{t}}
\een
for some polynomials
$q(z)\in\CC[z]$ and positive integers $t$ such that $t-\deg q\geqslant 2$ and $t\leqslant r+1$
so that $g_r(u)\in \CC((u))$ does exist and takes the required form.

Thus, by taking $c_0=1$ in the above argument, we may conclude that equation \eqref{g1} has a
solution $g(u)\in 1+(h/u)\tss\CC[[h/u,u]]$.
As shown in the proof of Proposition~\ref{prop:normal}
in Section~\ref{sec:QUE}, by multiplying $g(u)$ by an
appropriate element of $\CC[[h]]$ we get
a solution of both \eqref{g6} and \eqref{g1}.
To show that
any such solution is determined uniquely, expand $g(u)$ as in \eqref{guexp}.
We find from
\eqref{g6} that
\beql{g7}
g_{0}(u)=1\Fand\sum_{k=0}^l g_{k}(u)\tss g_{l-k}(-u)=0
\quad\text{for all}\quad l\geqslant 1.
\eeq
Returning to the induction argument in the first part of the above proof,
assume that the coefficients $g_{0}(u),\dots,g_{r-1}(u)\in\CC((u))$
are uniquely determined for some $r\geqslant 1$.
Relation \eqref{g3} with $l=r+1$
now determines $g_r(u)$ uniquely, up to an additive constant. However, its value is
fixed by the second condition in \eqref{g7} for $l=r$.


\bigskip\bigskip

\small


\noindent
S.K. \& A.M.:\newline
School of Mathematics and Statistics\newline
University of Sydney,
NSW 2006, Australia\newline
kslaven@maths.usyd.edu.au\newline
alexander.molev@sydney.edu.au\\

\noindent
S.K.:\\
Department of Mathematics\\
University of Zagreb, 10000 Zagreb, Croatia

\end{document}